\documentstyle[11pt]{article}
\font\eufm=eufm10
\def\frak#1{\hbox{\eufm#1}}
\newcommand{\bd}{
\begin{document}}
\newcommand{\ed}{\end{document}}
\newcommand{\be}{\begin{enumerate}}
\newcommand{\ee}{\end{enumerate}}
\newcommand{\bi}{\begin{itemize}}
\newcommand{\ei}{\end{itemize}}
\newcommand{\ba}{\begin{array}}
\newcommand{\ea}{\end{array}}
\newcommand{\sd}{\hspace{0.3ex}\tiny{\rhd\mbox{\hspace{-2ex}}<}\hspace{0.3ex}}
\newcommand{\g}{\frak g}
\newcommand{\ab}{\frak a}
\newcommand{\bb}{\frak b}
\newcommand{\h}{\frak h}
\newcommand{\k}{\frak k}
\newcommand{\al}{\frak a}
\newcommand{\n}{\frak n}
\newcommand{\p}{\frak p}
\newcommand{\q}{\frak q}
\newcommand{\m}{\frak m}
\newcommand{\ms}{\oplus}
\newcommand{\mt}{\otimes}
\newcommand{\dz}{\wedge}
\newcommand{\lra}{\longrightarrow}
\newcommand{\la}[2]{\Lambda_{#1#2}}
\newcommand{\veee}{{\tiny\mbox{$\vee$}}}
\newcommand{\kad}{ad^{\veee}}
\newcommand{\tf}{\tilde{f}}
\newcommand{\tu}{\tilde{u}}
\newcommand{\tg}{\tilde{g}}
\newcommand{\ove}{\overline}
\newtheorem{defi}{Definition}
\newtheorem{tw}{Theorem}
\newtheorem{prop}{Proposition}
\newtheorem{re}{Remark}
\newtheorem{lem}{Lemma}
\renewcommand{\arraystretch}{0.8}
\title{Double Lie algebras and Manin triples.}
\author{Piotr Stachura\thanks{Department of Mathematical Methods in Physics,
University of Warsaw, Ho\.{z}a 74, 00-682 Warszawa, Poland}}

\bd

\maketitle
\noindent{\bf Abstract\hspace{1em}} {\small The construction of {\em Lie
bialgebra} from {\em double Lie algebra} is presented. It is used to relate
some types of cobrackets on inhomogenous $so(p,q)$ algebras with double Lie
algebra structures on $so(p+1,q)$ and $so(p,q+1)$. Also is shown that the 
cobracket corresponding to $\kappa-deformation$, gives rise to the complete
Poisson-Lie Euclidean groups and non-complete Poincare groups.}

\section{Introduction}
Let $(G;A,B)$ be a double Lie group \cite{WeLu} i.e. $G$ is a Lie group, $A,B$
are closed subgroups and any element of $G$ has a unique decomposition: 
$g=a b\,,\,a\in A\,,\,b\in B$. Any double Lie group leads to a Manin group
and hence a pair of Poisson-Lie groups in duality (we do not require that
$G^*$ is simply connected). Let us recall that a Manin
group $(M;P,Q)$ is a double Lie group, where $M$ is equipped with invariant,
non-degenerate scalar product, vanishing on $TP$ and $TQ$.

We briefly sketch the way from double Lie groups to Manin groups\cite{SZ3}. 
Having a double Lie group we can define two compatible differential groupoid
structures on $G$ with $A$ and $B$ as sets of identities. This forms a 
$D^*$-group (in terminology of \cite{SZ3}). Applying the phase functor we get
two compatible symplectic groupoid structure on $T^*G$ ($S^*$-group). Then
using the symplectic form one can define invariant, non degenerate, scalar
product on $T^*G$. The sets of identities for both structures (namely
$P:=(TA)^0,Q:=(TB)^0$ )  are then Poisson-Lie groups, dual to each other.  

The infinitesimal version of a double Lie group is a double Lie algebra and
that of a Manin group is a Manin triple (cf. section 2). It is therefore
clear, that there should be a procedure which assigns a Manin triple (or a Lie
bialgebra) to each double Lie algebra. We present it in  section 2. 

In section 3 we relate some decompositions of $so(p+1,q)$ or $so(p,q+1)$
algebras with series of Lie bialgebra structures on $iso(p,q)$, among them
the so-called $\kappa$-deformation.

The main application of this study is to make a step towards a construction
on the $C^*$-algebra level of quantum groups corresponding to those Lie
bialgebras. By looking which double
Lie algebras give rise to (global) double Lie groups, we are in position to
distinguish between "good"(complete) and "bad"(non-complete)
cases \cite{SZ4}. In section 4 we show examples of decompositions from section
3 which do not give rise to a global decomposition of the corresponding Lie
groups. This is the case of the $\kappa$-deformation of Poincare group. It
strongly suggests that the $\kappa$-deformed Poincare group does not exist on
the $C^*$-algebra level. Contrary to this, the case corresponding to
$\kappa$-deformation of the Euclidean group is a "good" one: the global
decomposition is just the Iwasawa decomposition. Passing from groupoids to
their $C^*$ algebras \cite{Con}\cite{Ren} we get $\kappa$-deformed Euclidean
group on the $C^*$-algebra level. This will be described
elsewhere\cite{PS}.\vspace{1ex} 

{\em Throughout this paper all vector spaces, Lie algebras, Lie
groups are real and $\oplus$ means (if used without any comment) direct sum
of vector spaces. }

\section{Double Lie algebras and Manin triples.}

Let $(G,\pi)$ be a Poisson-Lie group with Lie algebra $\g$. Then $\pi$
determines linear mapping $\delta\,:\,\g\lra \g\dz\g $ which is a 1-cocycle
on $\g$ relative to adjoint representation on $\g\dz\g$ and the dual map
$\g^*\dz\g^*\lra \g^* $ is a Lie bracket. Conversely, if $G$ is connected and
simply connected then any such $\delta$ gives us a multiplicative Poisson
structure on $G$. Let us recall the following:
\begin{defi}
{\em \cite{WeLu} A pair $(\g,\delta)$ is said to be a }Lie bialgebra {\em if
$\g$ is a Lie algebra, $\delta\,:\,\g\lra\g\dz\g$ is a 1-cocycle on $\g$
relative to the adjoint representation of $\g$ on $\g\dz\g$ and
$\delta^*\,:\,\g^*\dz\g^*\lra\g^*$ is a Lie bracket.}
\end{defi}
In this situation we also say that $\delta$ is a cobracket on $\g$.
\begin{tw}
{\em (Manin)\cite{Dr}}
Let $\g$ be a Lie algebra, $\g^*$ its dual space and let $<\,,\,>$ denote the
canonical symmetric bilinear form on $\g\ms\g^*$. Let $\g^*$ be given a Lie
algebra structure. Then the dual map to the bracket on $\g^*$ is a cobracket
on $\g$ iff there exists a Lie algebra structure on $\g\ms\g^*$ such that:
\begin{enumerate}
\item $\g\,,\,\g^*$ are subalgebras of $\g\ms\g^*$.
\item The form $<,>$ on $\g\ms\g^*$ is invariant.
\end{enumerate}
In this case the bracket on $\g\ms\g^*$ is unique and is given by:\\
\mbox{$[X+\alpha,Y+\beta]=[X,Y]-\kad_{\beta} X +\kad_{\alpha}Y
+[\alpha,\beta]+\kad_X\beta-\kad_Y\alpha$} where $[\,,\,]$ denotes brackets
on $\g$ and $\g^*$ and $\kad$ denotes the coadjoint representations of $\g$
and $\g^*$. The $\g\ms\g^*$ with the above bracket will be denoted by
$\g\bowtie \g^*$.
\end{tw}
The Lie bialgebras are in one to one correspondence with Manin triples.
\begin{defi}
{\em \cite{Dr} Three Lie algebras $(\m;\p,\q)$ form a} Manin triple {\em
iff: \begin{enumerate}
\item $\p,\q$ are Lie subalgebras of $\m$ and $\m=\p\ms\q$ as a vector
space sum. 
\item $\m$ is equipped with invariant, non-degenerate scalar product such
that $\p,\q$ are isotropic. 
\end{enumerate}}
\end{defi}
We need also the definition of a double Lie algebra:
\begin{defi}{\em \cite{WeLu}}
A double Lie algebra {\em is a triple $(\g;\ab,\bb)$ such that $\ab,\bb$ are
Lie subalgebras of $\g$ and $\g=\ab\ms\bb$ as a vector space sum.}
\end{defi}

Now let $(\g;\ab,\bb)$ be a double Lie algebra. Consider the coadjoint action
semi-direct product $\g\sd\g^*$ with natural bilinear, symmetric, 
invariant form $<\,,\,>$. 
Then $\ab^0$ is $\ab$ invariant and $\bb^0$ is $\bb$ invariant (where
$\ab^0, \bb^0$ denotes annihilators of $\ab$ and $\bb$ respectively). \\
Indeed if \mbox{$x,y\in\ab$}, \mbox{$\alpha\in\ab^0$} then
\mbox{$<[x,\alpha],y>$}=\mbox{$<\alpha,[y,x]>=0$} 
since $\ab$ is a subalgebra. Of course $\ab\ms\ab^0$ and $\bb\ms\bb^0$ are
isotropic. In this way we get that $(\g\sd\g^*; \ab\sd\ab^0,\bb\sd\bb^0)$ is
a Manin triple and we have a bialgebra structure on $ \ab\sd\ab^0.$ 
The cobraket $\delta$ satisfies: $\delta(\ab)\subset \ab^0\dz\ab$ and
$\delta(\ab^0)\subset\ab^0\dz\ab^0$. So \mbox{$(\ab\sd\ab^0)\bowtie
(\ab\sd\ab^0)^*= \g\sd\g^*$} where we identified $\bb\ms\bb^0$ with
$(\ab\ms\ab^0)^*$. Thus we have a procedure which associates with each double
Lie algebra a Lie bialgebra $(\h,\delta)$ with the following properties:
1. $\h=\ab\sd V$ ($V$-abelian ideal), 2. $\delta(\ab)\subset\ab\dz V,
\delta(V)\subset V\dz\ V.$

Conversely, let $(\h:=\ab\sd V,\delta)$ be a semidirect product with abelian
ideal $V$ and cobraket $\delta$ such that $\delta(\ab)\subset V\dz \ab$ and
$\delta(V)\subset V\dz V$. $\delta$ defines a Lie algebra structure on
$\h^*=\ab^0\ms V^0$. Let us show, that $\h^*=\ab^0\sd V^0$ with abelian ideal
$V^0$. We adopt the following notation: capital letters $A,B,...$ are elements
of $\ab$, $A^0,B^0,...$ are elements of $\ab^0$, small $x,y,...$ are elements of
$V$ and $x^0,y^0,...$ are elements of $V^0$. We have: 
\bi
\item $<C,[A^0,B^0]>=<\delta(C),A^0\dz B^0 >=0$, since $\delta(\ab)\subset V\dz
\ab$, so  $\ab^0$ is a subalgebra; 
\item $<z,[x^0,y^0]>=<\delta(z),x^0\dz y^0 >=0 $, since $\delta(V)\subset V\dz
V \\
<C,[x^0,y^0]>=<\delta(C),x^0\dz y^0>=0$, so $V^0$ is an abelian subalgebra;
\item  $<z,[A^0,x^0]>=<\delta(z),A^0\dz x^0>=0$ and $V^0$ is an ideal.
\ei
Now we prove that the Lie algebra $\h\ms\h^*$ with bracket given in the
theorem 1 coincides with semidirect product $(\ab\ms \ab^0)\sd(\ab\ms \ab^0)^*$
with coadjoint action if we identify $V\ms V^0$ with $(\ab\ms \ab^0)^*$ by
duality: $<<A+B^0,x+y^0>>:=<A,y^0>+<x,B^0>$.  \\
For the coadjoint representation of $\h$ we have:
\bi
\item $\kad_{\ab}(V^0)\subset V^0$ : $<z,\kad_A(y^0)>=<[z,A],y^0>=0$ (since
$[\ab,V]\subset V$);
\item $\kad_{\ab}(\ab^0)\subset \ab^0 $ : $<C,\kad_A(B^0)>=<[C,A],B^0>=0$ (since
$[\ab,\ab]\subset \ab$ );
\item $\kad_V(V^0)=0 $ : $<z,\kad_x(y^0)>=<[z,x],y^0>=<0,y^0>=0\,,\\
<A,\kad_x(y^0)>=<[A,x],y^0>=0$ (since $[\ab,V]\subset V$);
\item $\kad_V(\ab^0)\subset V^0$ : $<z,\kad_x(A^0)>=<[z,x],A^0>=<A^0,0>=0$.
\ei
And for the coadjoint representation of $\h^*$:
\bi
\item $\kad_{\ab^0}(\ab)\subset \ab$ : $<\kad_{A^0}(B),C^0>=<B,[C^0,A^0]>=0$ (since
$[\ab^0,\ab^0]\subset \ab^0$);
\item $\kad_{\ab^0}(V)\subset V$ : $<\kad_{A^0}(y), z^0>=<y,[z^0,A^0]>=0$
(since $[\ab^0,V^0]\subset V^0)$;
\item $\kad_{V^0}(\ab)\subset V$ : $<\kad_{y^0}(A),z^0>=<A,[z^0,y^0]>=0$ (since
$[V^0,V^0]=0$);
\item $\kad_{V^0}(V)=0$ : $<\kad_{x^0}(y),z^0>=<y,[z^0,x^0]>=<y,0>=0\,,\\
<\kad_{x^0}(y),A^0>=<y,[A^0,x^0]>=0$ (since $[\ab^0,V^0]\subset V^0$).
\ei
From this it follows that $V \ms V^0$ is an abelian ideal and $\ab\ms \ab^0$ is a
subalgebra. If we identify $V\ms V^0$ with $(\ab\ms \ab^0)^*$ we are in the
situation in Theorem 1. and it follows that the action of $\ab\ms \ab^0$ is a
coadjoint action. In this way we have proved the following:
\begin{prop}
Any double Lie algebra $(\g;\ab,\bb)$ leads to a Lie bialgebra $(\h,\delta)$
such that $\h=\ab\sd V$ is a semi-direct product with abelian ideal $V$ and
the cobracket $\delta$ satisfies: $\delta(\ab)\subset \ab\dz V$ and
$\delta(V)\subset V\dz V$. Conversely, any Lie bialgebra of this type is
obtained in this way.  
\end{prop}

\section{Iwasawa-type decompositions of $so(p,q)$  and
bialgebra structures on $iso(p-1,q)\,,\,iso(p,q-1)$.}
\subsection{Inhomogenous $so(p,q)$ algebras and $b$-type Poisson structures.}

Let $(V,\eta)$ be a n+1 dimensional, real vector space with symmetric,
nondegenerate, billinear form $\eta$ of signature $(p,q)$. By $\eta$ we also
denote isomorphism $V \lra V^*$ given by $\eta(x)(y):=\eta(x,y)$.  
Let $iso(p,q):=so(p,q)\sd V$ be an inhomogenous $so(p,q)$ Lie algebra.
$so(p,q)$ is isomorphic to $V\dz V$ by: $x\dz
y\mapsto \la{x}{y}:=x\mt\eta(y) -y\mt\eta(x)$. If $(e_i)$ is an orthonormal
basis  of $V$: $(\la{i}{j}:=\la{e_i}{e_j}\,,\,i<j) $ form a basis of
$so(p,q)$, with commutators:  $[\la{i}{j},\la{k}{l}]=\eta_{il}
\la{j}{k}+\eta_{jk}\la{i}{l}-\eta_{ik}\la{j}{l}- \eta_{jl}\la{i}{k}$.\\
Let $K(A,B):=-\frac{1}{2} Tr(A B)$. This is $ad$ invariant, non degenerate
scalar product on $so(p,q)$ and $(\la{i}{j}\,,\,i<j)$ form an orthonormal
basis. \\
Let $(\la{i}{j}^*\,,\,e_k^*\,:\,i <j)$ be a
basis in $iso(p,q)^*$ given by:  \\
$<\la{k}{l},\la{i}{j}^*>:=K(\la{i}{j},\la{k}{l})=\eta_{ik}\eta_{jl}-\eta_{il}
\eta_{jk}$ and $<e_l,e_k^*>:=\eta_{kl}$. \\
In this basis the coadjoint representation of $iso(p,q)$ has the following
form: \\
$\kad_{\la{a}{b}}(\la{c}{d}^*)=\eta_{ad}\la{b}{c}^*+\eta_{bc}\la{a}{d}^*-
\eta_{ac}\la{b}{d}^*-\eta_{bd}\la{a}{c}^*\;,\;
\;\kad_{\la{a}{b}}(e_k^*)=\eta_{kb} e_a^* - \eta_{ka} e_b^*\,,\\
\kad_{e_a}(\la{c}{d}^*)=0\;,\;\; \kad_{e_a}(e_b^*)=-\la{a}{b}^*. $

Let $\g:=iso(p,q)$ and $\h:=so(p,q)$, so $\g=\h\sd V$.
It is known \cite{SZ1} that all bialgebra structures on $\g$ for $p+q>2$
are coboundary i.e. are of the form $\delta=\partial r$ for some $r\in
\g\dz\g$ ($\delta(x)=\partial r (x):=ad_x(r)$) where $r$ satisfies the
generalized classical Yang-Baxter equation: $[r,r]\in (\g\dz\g\dz\g)_{inv}$.
Since $\g\dz\g=(\h\dz\h)\ms(\h\dz V)\ms (V\dz V)$ we can write :
$r=c+b+a\,,\, c\in \h\dz\h\,,\,b\in\h\dz V\,,\,a\in V\dz V$.  We say that $r$
is of $b$-type iff $r=b$. In this case $b$ satisfies $[b,b]=t\Omega\,,\,t\in
R$ where $\Omega:= \eta^{jl}\eta^{km}e_j\dz e_k\mt\la{l}{m}$ is the canonical
$\g$-invariant element of $\g\dz\g\dz\g$.

We will be interested in the following solutions of this equation:
\be
\item $b_x:=\eta^{jk}e_j\dz\la{x}{e_k}\,,\,x\in V\,$ 
is a solution  with $t=-\eta(x,x)$ \cite{SZ1}.
\item $\tilde{b}_x:=b_x+x\dz X$ where $X \in\h$ and $X x=0$ is a solution
with the same $t$ \cite{SZ1}. 
\item Let $x\in V$ be a null vector and let $v_i\in V,\,X_i\in \h$ satisfy:
$X_i x=0,\,X_i v_j=-\delta_{ij} x,\,[X_i,X_j]=0$. Then $b:=b_x+x\dz Y+\sum
v_i\dz X_i$, where $Y:=\sum \alpha_i X_i,\,\alpha_i\in R$ is a solution
with $t=0$ \cite{PS1}. 
\item $b:=\tilde{b}_x+v\dz X$ where $Xv=v$ is a solution with $t=-\eta(x,x)$ \cite{PS1}.
\end{enumerate}

We will need the fact that $b$ is completely determined by
the bracket on $V^*$\cite{SZ2}. 
Let $b=v_i\dz h_i$ (sumation). We use the same letter for the mapping $b: V^*
\lra \h $ given by: $ b(\alpha):=<v_i,\alpha>h_i$. Then the bracket on $V^*$
can be expressed by $b:\,[\alpha,\beta]=b(\alpha)\beta-b(\beta)(\alpha)$,
where the action of $\h$ on $V^*$ is a coadjoint action.
Let $e_k^*:=\eta(e_k)$ and $b(e_k^*)=:{b_k}^{mn}\la{m}{n}$ with
${b_k}^{mn}=-{b_k}^{nm}$, $[e_i^*,e_j^*]=:{f_{ij}}^k
e_k^*,\,{f_{ij}}^k=-{f_{ji}}^k$. Then: 
$[e_i^*,e_j^*]={b_i}^{kl}(\eta_{lj}e_k^*-\eta_{kj}e_l^*)-{b_j}^{kl}
(\eta_{il}e_k^*- \eta_{ik}e_l^*)=\\=2({b_i}^{kl}\eta_{lj}-{b_j}^{kl}\eta_{il})
e_k^*={f_{ij}}^k e_k^*$. From this it follows that
$f_{ijk}=2(b_{jik}-b_{ijk})$ (we used $\eta_{ij}$ to lower indeces).
This equation determines $b_{ijk}$: \mbox{$b_{ijk}=\frac{1}{4}
(f_{jki}-f_{ijk}-f_{kij})$} and $b$, since $b=b^{kmn}e_k\dz\la{m}{n}$.\\
Let us also notice the following:
\begin{lem} Let $p+q>2$ and let $\delta=\partial r$ be a cobracket on
$iso(p,q)=so(p,q)\sd V=:\h\sd V$ which satisfies: $\delta(\h)\subset \h\dz V
,\,\delta(V)\subset V\dz V$. Then $r=b$.  \\
Proof: {\em $r=c+b+a\,,\, c\in \h\dz\h\,,\,b\in\h\dz V\,,\,a\in V\dz V$. Then
the condition  $\delta(\h)\subset \h\dz V$ is equivalent to
$ad_h(c)=0$ and $ad_h(a)=0$ for $h\in\h$ and $\delta(V)\subset V\dz V$ is
equivalent to $ad_v(c)=0$ for $v\in V$. But since $\h$ is semisimple, from the
first equality it follows that $c=0$. Also since isomorphism $V\dz V\ni x\dz
y\mapsto\la{x}{y}\in\h$ intertwines action of $\h$ on $V\dz V$ with the
adjoint action on $\h$ we have $a=0$. $\Box$}
\end{lem}

\subsection{Iwasawa-type decomposition of $so(p,q)$.}
In this section we put: $p+q=:n+1\,,\;n>2;$\\
$so(p,q)\supset\h_1:=<\la{i}{j}\,:\,2\leq i < j\leq n+1>=so(p-1,q);\\
so(p,q)\supset\h_2:=<\la{i}{j}\,:\,1\leq i < j\leq n>=so(p,q-1)$;\\
$f:=\la{1}{n+1}\;,\;\;\;g_k:=\la{1}{k}+\la{k}{n+1}\,,\;2\leq k\leq n
\;,\;\n:=<g_2,...,g_n>$.\\ 
With this notation: $so(p,q)=\h_1\ms<f>\ms \n=\h_2\ms<f>\ms \n$.\\
Then $u:=<f>\ms \n$ is a Lie subalgebra:
$[f,g_k]=\eta_{11}g_k=g_k\;,\;[g_k,g_l]=0$.  
The corresponding annihilators in $so(p,q)^*$ are equal:
$\h_1^0=<\la{1}{l}^*\,,\,2\leq l\leq n+1>\,,
\,\,\h_2^0=<\la{l}{n+1}^*\,,\,1\leq l\leq n>\,,\,\,\\
u^0=<g_l^*\,,\,2\leq l\leq n>\ms <\la{m}{s}^*\,,\,2\leq m,s\leq n>$ where
$g_l^*:=\la{1}{l}^*+\la{l}{n+1}^*$. \\
In this way we get two double Lie algebras: $(so(p,q);so(p-1,q),u)$ and
$(so(p,q);so(p,q-1),u)$.\vspace{1ex}

{\it The double Lie algebra $(so(p,q);so(p-1,q),u)$.}\\
We use $K$ to equip $\h_1^0$ with scalar product. With respect to this
product $(\la{1}{l}^*\,,\,2\leq l\leq n+1)$ form an orthonormal basis. The
signature is $(p-1,q)$. 
The action of $\h_1$ is given by:\\
$\la{i}{j}(\la{1}{l}^*)=\eta_{jl}\la{1}{i}^* -\eta_{il}\la{1}{j}^*$, so
$\h_1\sd\h_1^0 =iso(p-1,q)$. Let us put $v_l:=\la{1}{l}^*$.  We already know
that $\delta=\partial b$ for some $b\in \h_1^0 \dz \h_1$ and $b$ is determined
by the bracket on $u=(\h_1^0)^*$. We have: $f=v_{n+1}^*\,,\,g_l=v_l^*$ and
$[v_l^*,v_m^*]=\delta_{ln+1} v_m^*-\delta_{mn+1} v_l^*$, so  
${f_{lm}}^s=\delta_{ln+1}\delta_m^s-\delta_{mn+1}\delta_l^s$.  Using
formula from section 2: 
$b_{ijk}=\frac{1}{2}(\delta_{jn+1}\eta_{ik}-\delta_{kn+1}\eta_{ij})$ 
and $b=\eta^{sl}v_s\dz\la{n+1}{l}=b_{v_{n+1}}$.  

In this way we have shown that the double Lie algebra $(so(p,q);so(p-1,q),u)$
leads to the cobracket on $iso(p-1,q)$ given by $b=b_x\,,\,\,\eta(x,x)<0$.
\vspace{1ex}

{\it The double Lie algebra $(so(p,q);so(p,q-1),u)$.}\\
Now we use $-K$ as a scalar product on $\h_2^0$. With respect to this product
product \mbox{$(\la{l}{n+1}^*,\,1\leq l\leq n)$} form an orthonormal basis. The
signature is $(p,q-1)$. The action of 
$\h_2$ is given by:\\
$\la{i}{j}(\la{l}{n+1}^*)=\eta_{jl}\la{i}{n+1}^* -\eta_{il}\la{j}{n+1}^*$, so
$\h_2\sd\h_2^0 =iso(p,q-1)$. Let us put $v_l:=\la{l}{n+1}^*$ then
\mbox{$f=-v_{1}^*$}, \mbox{$g_l=-v_l^*$.} 
In the same way as above we get:
${f_{lm}}^s=\delta_{1m}\delta_l^s-\delta_{1l}\delta_m^s\,,\,\,
b_{ijk}=\frac{1}{2}(\delta_{1k}\eta_{ij}-\delta_{1j}\eta_{ik})$ and 
$b=-\eta^{sm}v_s\dz\la{1}{m}=-b_{v_1}$. 

So we see that the double Lie algebra $(so(p,q);so(p,q-1),u)$
leads to the cobracket on \mbox{$iso(p,q-1)$} given by $b=b_x\,,\,\,\eta(x,x)>0$.
\vspace{1ex}

We can be a little bit more general and instead of $f$ take
$\tf:=\la{1}{n+1}+s$ where\\ $s:=s^{ij}\la{i}{j} \in
<\la{i}{j}\,,\,2\leq i,j\leq n>=so(p-1,q-1)\,,\,s^{ij}=-s^{ji}$.
Then $[\tf,g_k]=g_k+s^{ij}(\eta_{jk}g_i- \eta_{ik}g_j)$ so again
$\tu:=<\tf>\ms \n$ is a Lie subalgebra and $so(p,q)=\h_1\ms\tu=\h_2\ms\tu$. 
Let us analyze the new cobrackets on $\h_1\sd \h_1^0$ and $\h_2\sd\h_2^0$.
\vspace{1ex}

{\it Double Lie algebra $(so(p,q);so(p-1,q),\tu)$.}\\
Keeping the same notation as above: $v_k^*=g_k\,,\,v_{n+1}^*=\tf$.\\
\mbox{$[v_l^*,v_m^*]=[\delta_{ln+1}(\delta_m^s+s^{sj}\eta_{jm}-
s^{is}\eta_{im})
-\delta_{mn+1}(\delta_l^s+s^{sj}\eta_{jl}-s^{is}\eta_{il})]v_s^*$}. 
From this it follows that:\\ \mbox{$b^{sml}=\frac{1}{2}(\eta^{mn+1}\eta^{sl}-
\eta^{ln+1}\eta^{sm})-\frac{1}{2}\eta^{n+1s}(s^{lm}-s^{ml})$} and\\
$b=\eta^{mp}v_m\dz \la{n+1}{p}+v_{n+1}\dz
s=b_{v_{n+1}}+v_{n+1}\dz s$. \\
So we get the cobracket on $iso(p-1,q)$ given by $b=b_x+x\dz X
\,,\,\,\eta(x,x)<0$.  
\vspace{1ex}

{\it Double Lie algebra $(so(p,q);so(p,q-1),\tu)$.}\\
Now
\mbox{$[v_l^*,v_m^*]=[\delta_{1m}(\delta_l^s+s^{sj}\eta_{jl}-s^{is}\eta_{il})
-\delta_{1l}(\delta_m^s+s^{sj}\eta_{jm}-s^{is}\eta_{im})]v_s^*$}. From this
equation: $b^{sml}=\frac{1}{2}(\eta^{1l}\eta^{sm}-
\eta^{1m}\eta^{sl})-\frac{1}{2}\eta^{1s}(s^{lm}-s^{ml}) $ and
$b=-b_{v_1}-v_1\dz s$.\\
This is the cobracket on $iso(p,q-1)$ given by $b=b_x+x\dz X
\,,\,\,\eta(x,x)>0$.  

In this way we have shown that bialgebra structures on $iso(p,q)$ of type 1
and 2 for non null vectors come from the double Lie algebra structures on
$so(p+1,q)$ or $so(p,q+1)$.
\vspace{1ex}

Now is time for type 4. Let $\tf=f+s$ be as above and let us assume that $s$
has d-dimensional eigenspace with eigenvalue 1. Then this is null subspace
and one can choose an orthonormal basis \mbox{$(e_i\,,\,2\leq i\leq n)$}
such that this eigenspace is equal:\\
\mbox{$<e_{m_1}-e_{n_1},...,e_{m_d}-e_{n_d}>$} with $2~\leq~m_k~<~n_k~\leq
~n\,,\,\eta_{m_km_k}=1\,,\,\eta_{n_kn_k}=-1$ for $k=1,...,d$. Let
$D:=\{m_1,n_1,...,m_d,n_d\}$ and $\chi^D$ be the characteristic function of
$D$. A short calculation shows that in this situation:
$s^{ij}\eta_{jp}\chi^D(i)-s^{ij}\eta_{ip}\chi^D(j)=-\chi^D(p)$.  \\
Let us define: $\tg_k:=\chi^D(k) f+g_k\,,\,k=2,...,n$ and
$\tilde{U}:=<\tf>\ms <\tg_2,...,\tg_n>.$ 
\begin{lem}
$so(p,q)=\h_1\ms\tilde{U}=\h_2\ms\tilde{U}$ and 
$\,\tilde{U}$ is a Lie subalgebra.\\
Proof: 
$[s,\tg_p]=[s,g_p]=s^{ij}(\eta_{jp}g_i-\eta_{ip}g_j)=
s^{ij}(\eta_{jp}(g_i+\chi^D(i) f)-\eta_{ip}(g_j+\chi^D(j)
f))+\\-(s^{ij}\eta_{jp}\chi^D(i) -s^{ij}\eta_{ip}\chi^D(j))f=
s^{ij}(\eta_{jp}\tg_i-\eta_{ip}\tg_j)+\chi^D(p)f.$\\
{\em From this it follows:} $[\tf,\tg_p]=[f+s,\chi^D(p) f+g_p]=
\tg_p+s^{ij}(\eta_{jp}\tg_i-\eta_{ip}\tg_j)$.\\ 
$[\tg_k,\tg_l]=\chi^D(k) g_l-\chi^D(l)g_k=\chi^D(k) \tg_l-\chi^D(l)\tg_k$\\
{\em In this way $\tilde{U}$ is a subalgebra and simple
calculations show that it is complementary to $\h_1$ and $\h_2.\; \Box$}
\end{lem}

{\it Double Lie algebra $(so(p,q);so(p-1,q),\tilde{U})$.}\\
We have: $\tf=v_{n+1}^*\,,\,\tg_l=g_l=v_l^*$ for $l\not\in D$ and 
$\tg_l-\tf=v_l^*$ for $l\in D$.
From the bracket on $\tilde{U}$: $[v_{n+1}^*,v_k^*]=v_k^*+
s^{ij}(\eta_{jk}v_i^*-\eta_{ik}v_j^*)$ and  
$[v_k^*,v_l^*]=2(\chi^D(l)\eta_{ki}s^{mi}-\chi^D(k)\eta_{li}s^{mi})v_m^*$.\\
So $f_{ijk}=F_{ijk}+2(\chi^D(j)s_{ki}-\chi^D(i)s_{kj})$ where we put
$F_{ijk}$ - the structure constants for\\ \mbox{$(so(p,q);so(p-1,q),\tu)$}.
In this way: $b=b_{v_{n+1}}+v_{n+1}\dz s-
((v_{m_1}-v_{n_1})+\dots +(v_{m_d}-v_{n_d}))\dz s$ and this is $b$ of type 4
for $\eta(x,x)<0.$
\vspace{1ex}

{\it Double Lie algebra $(so(p,q);so(p,q-1),\tilde{U})$.}\\
This is completely analogous to the above case and we get solution of type 4
with $\eta(x,x)>0$.
\begin{re}{\em The solutions of type 1 and 2 for null vectors and of type 3
are the special cases of the following double Lie algebras. \\
Let $\ove{f}:=e_1+\lambda \la{1}{n+1} +s +g$ where $\lambda\in R\,,\,s$-as
above$,\, g\in<\la{1}{k}+\la{k}{n+1}\,,\,k=2,...,n>$ \\
$w:=e_1-e_{n+1}.$ For any basis $(x_k)$ of $<e_2,...,e_n>$ let
$g_k:=\la{1}{x_k}+\la{x_k}{n+1}.$ Then $<\ove{f},w,x_k>$ is a subalgebra of
$iso(p,q)$ complementary to $so(p,q)$, the same holds for
$<\ove{f},w,x_k+\lambda g_k>.$ Moreover if $<e_2,...,e_n>$ is a direct sum
of $s$-invariant subspaces we can make above choice on each subspace
separately. This leads to the family of double Lie algebras of form
$(iso(p,q);so(p,q),\ab)$ and a family of cobrackets on $so(p,q)\sd
so(p,q)^0\subset iso(p,q)\sd (iso(p,q))^*$ which we can identify with
$iso(p,q)$. }\end{re}
\section{Global decompositions}
Let $G$ be a Poisson-Lie group. Then $\g^*$ and $\m:=\g\bowtie\g^*$ are Lie
algebras. We consider the following problem: to find a connected Lie group
$M$ with Lie algebra $\m$ such that:
\begin{enumerate}
\item $G$ is a Lie subgroup of $M$
\item $M=G G^*$ (or at least $G G^*$ is dense in $M$), where $G^*$ is the
analytic subgroup  of $M$ with Lie algebra $\g^*$.
\end{enumerate}

We study this problem for two of the Poisson-Lie groups obtained in the last
section: $E(n):=SO(n)\sd R^n$ and $P_0(n):=SO_0(1,n-1)\sd R^n$ coming from the
double Lie algebras $(so(1,n);so(n),\tu)$ and $(so(1,n);so(1,n-1),\tu)$
(notation as in section 3).

Let $V$ be a finite dimensional, real vector space and $K\subset GL(V)$ a
closed, connected subgroup which acts on V without fixed points (except 0).
Let $G:=K\sd V$ be a semidirect product and $\g=\k\sd v$ its Lie algebra. In
this situation the center of $G$ is trivial and $G=Int(\g)$ ($Int(\g)$ is the
adjoint group of $\g$\cite{He}.) 

Suppose we are given a bialgebra structure on $\g$ as in 
section 2. Then we know that $\g^*=\k^0\sd v^0$ with $v^0$-abelian ideal 
and $\m:=\g\bowtie \g^*=(\k\ms \k^0)\sd (v\ms v^0)=:\h\sd \h^*$ and the
action is a coadjoint action. We assume that $\h$ is semisimple. 
Then the center of $\m$ is trivial. 

Let $M$ be a connected Lie group with Lie 
algebra $\m$; $H,H^*,G^*,K^0,V^0$ be  analytic subgroups with Lie algebras 
$\h,\h^*,\g^*,\k^0,v^0$ respectively. Moreover let us assume that $G$ is
contained in $M$, so $G,K,V$ are identified with analytic subgroups of $M$
with Lie algebras $\g,\k,v$ and $K,V$ are closed in $G$.  Since $\h$ is 
semisimple $\tilde{H}:=Int(\h)$ is a closed \cite{He} (and by definition
connected) 
subgroup of $GL(\h)$ and $\tilde{H}$ acts on $\h^*$ without fixed points
(except 0). Let $\tilde{M}:=\tilde{H}\sd \h^*$ then $\tilde{M}$ has
trivial center and the same Lie algebra as $M$, so $\tilde{M}=Int(\m)$ and
$M$ is a covering group of $\tilde{M}$, the covering homomorphism is given by 
$\phi:=Ad_M$. We have $\tilde{H}=\phi(H)$ and let
$\tilde{G^*},\tilde{K^0},...$ denote the images by $\phi$ of $G^*,K^0,...$.  
These are analytic subgroups of $\tilde{M}$ with Lie algebras
$\g^*,\k^0,...$. We are going to prove the following:
\begin{prop} If $\,G G^*$ is dense in $M$ then $\tilde{K} \tilde{K^0}$ is
dense in $\tilde{H}$.

Proof: {\em Let $\tilde{h}\in\tilde{H}$ and $\tilde{h}=\phi(m)$ for some
$m\in M$. 
$G G^*$ is dense in $M$, so $m=\lim g_n g_n^*$ for $g_n\in G,\,g_n^*\in
G^*$.  But $g_n=k_n v_n,\,k_n\in K,\,v_n\in V$ and since $\g^*=\k^0\sd v^0$
any element of $G^*$ has, possibly non unique, decomposition $g_n^*=k_n^0
v_n^0,\,k_n^0\in K^0,\,v_n^0\in V^0$. Because $\h^*$
is an ideal in $\m$, $H^*$ is normal subgroup and $v_nk_n^0=k_n^0 x_n$ for
some $x_n\in H^*$. In this way $\tilde{h}=\lim \phi(k_n)\phi(k_n^0)\phi(x_n)
\phi(v_n^0)=\lim \phi(k_n)\phi(k_n^0)\phi(x_n v_n^0)$ with
$\phi(k_n)\phi(k_n^0)\in \tilde{K}\tilde{K^0}\subset \tilde{H}$ and $\phi(x_n v_n^0) \in
\h^*$. Now, the convergence 
in $\tilde{M}$ is convergence along "coordinates" in $\tilde{H}$ and $\h^*$
so $\tilde{h}=\lim \phi(k_n) \phi(k_n^0)$ and $\tilde{K}\tilde{K^0}$ is dense
in $\tilde{H}$. $\Box$ }
\end{prop}
In the following we need the {\it Iwasawa decomposition of
$\,SO_0(1,n)$.} \cite{He} \\
Let $so(1,n)=\k\ms\al\ms\n$ where $\k:=<\la{i}{j}\,,\,2\leq i,j\leq
n+1>=so(n) \,,\; \al:=<\la{1}{n+1}>\,,\;\\
\n:=<\la{1}{k}+\la{k}{n+1}\,,\,2\leq
k\leq n>$ be the Iwasawa decomposition of $so(1,n)$. To this corresponds
decomposition of a connected component of the identity: $SO_0(1,n)=K A N$
where $K,A,N$ are analytic subgroups of $SO_0(1,n)$ with algebras $\k,\al,\n$
respectively. \vspace{1ex}
\\In our case these subgroups look as follows: 
$K=\left\{\left( \ba{cc} 1 & 0\\0 & T \ea \right)\;,\,T\in
SO(n)\right\}$ \\
$A$ is one parameter subgroup: $A(t)=\left( \ba{ccc} \cosh t & 0
&\sinh t\\ 0 & I &0\\ \sinh t & 0 &\cosh t \ea \right) \,,\, t\in R\,,\, I$
is $(n-1)\times (n-1)$ identity matrix; \vspace{1ex}\\ 
and elements of $N$: $N(x)= \left( \ba{ccc}
1+\frac{1}{2} |x|^2 & -x & \frac{1}{2}|x|^2 \\ -x^t &
I & -x^t \\  
-\frac{1}{2}|x|^2 & x & 1-\frac{1}{2}|x|^2 
\ea \right)\,,\; x:=(x_2,...,x_n)\in
R^{n-1},\,|x|^2~=~\sum_{i=2}^n x_i^2$.\vspace{1ex}\\
Moreover $N$ is commutative
and $A(t) N(x)=N(e^{-t}x) A(t)$.\vspace{1ex}

Now we pass to  Poisson-Lie groups
$E(n):=SO(n)\sd R^n$ and $P_0(n):=SO_0(1,n-1)\sd R^n$. 
In both cases $\tilde{H}=SO_0(1,n)\,,\,
\tilde{M}=SO_0(1,n)\sd so(1,n)^*$ ;
$\k^0=\tu=<\tf>\ms\n\,,\\
\tilde{K^0}=F N$ where 
$F$ is one parameter subgroup of elements: $\exp (t \tf)=:F(t)= A(t) S(t)
=S(t) A(t)$,\\
$S(t):=\left(\ba{ccc} 1 & 0 & 0\\ 
0 &\exp (t s)&0\\0&0&1\ea \right)\;\;,\,
A(t),\,N$-as above.
\vspace{1ex}

{\it The Euclidean group.} \\
$\;\k:=so(n)=<\la{i}{j}\;,\,2\leq i,j\leq n+1>\,,\,\tilde{K}=\left\{\left(
\ba{cc} 1 & 0\\0 & T \ea \right)\;,\,T\in SO(n)\right\}$. In this case we will show that the
global decomposition $SO_0(1,n)=\tilde{K} F N$ holds and this is
sligthly modified Iwasawa decomposition.\\
Let $SO_0(1,n)\ni g=k A(p) N(y)\,,\,\,k\in K\,,\,A(p)\in
A\,,\,N(y)\in N$ be a decomposition of $g$. We seek for
$\tilde{k}\in\tilde{K}=K\,,\, F(t)\,,\,N(x)$ such that:
$\tilde{k} F(t) N(x)=g= k A(p) N(y)$. Using the definition of $F$
: $\tilde{k} S(t) A(t) N(x)= k A(p) N(y)$. Since $S(t)\in K$ and
the decomposition is unique we have: $\tilde{k} S(t)=k\,;\;
t=p\,;\,x=y$ and $\tilde{k}=k S(-p)$.  
This proves that the decomposition is global and if $s=0$ this is just the
Iwasawa decomposition. Hence the Manin group for $E(n)$ is $SO_0(1,n)\sd
so(1,n)^*$. 

{\it The Poincare group.}\\ $\;\k:=so(1,n-1)=<\la{i}{j}\;,\,1\leq i,j\leq n>
\,,\,\tilde{K}=\left\{\left( \ba{cc} \tilde{T} & 0\\0 & 1\ea 
\right) \;,\,\tilde{T}\in SO_0(1,n-1)\right\}$.\\ We will show that
$SO_0(1,n)\setminus \tilde{K} F N$ contains an open subset.\\
Let $W$ be a function on $SO_0(1,n)$ defined by:
$W(g):=\eta(g(e_1-e_{n+1}),e_{n+1})$. This function is obviously continous,
and is easy to see that if $g=\tilde{k} f n\,,\,\tilde{k}\in\tilde{K},f\in
F,n\in N$ then $W(g)>0$. 
But for $SO_0(1,n)\ni k_0:=\left(\begin{array}{cc} I_{n-1} &0\\
0 & -I_2\end{array}\right)\,,\,I_l$-is $l\times l$ identity matrix we have
$W(k_0)<0$ what proves the assertion.

Next we show that one can find non-connected extension of $P_0(n)$ for which
there exists connected (in fact simply connected) Poisson dual group $G^*$
and the set $G G^*$ is dense in $M$.

\begin{lem} $\tilde{K} F N=\{k a n\,:\,k\in K,a\in A,n\in
N,k_{n+1 n+1}>0\}.$\\ 
Proof: {\em We try to solve: 
$\tilde{k} F(t) N(x)= k A(p) N(y) \;,\,\tilde{k}\in
\tilde{K}\,,\,k\in K$. 
Using the commutation relation between $A$ and $N$ we
get:\\ $\tilde{k} S(t)=k A(p) N(y-x)A(-t)=k A(p-t)
N(e^{-t}(y-x))$.\\ Let $z:=e^{-t}(y-x)$ and
$w:=p-t$, then:\vspace{1ex}\\
 $A(w) N(z)=\left( \ba{ccc} \cosh w+\frac{1}{2}|z|^2
e^{-w} & -z e^{-w} & \sinh w+\frac{1}{2} |z|^2 e^{-w}\\
-z^t & I & -z^t\\ \sinh w-\frac{1}{2}|z|^2 e^{-w} &z
e^{-w} & \cosh w-\frac{1}{2}|z|^2 e^{-w} \ea \right).$\vspace{1ex}\\
Now we look at the last row of the equality $\tilde{k} S(t)=k A(w) N(z)$
: \\
$(n+1,1):\,-\sum_{j=2}^n k_{n+1j}z_j+k_{n+1n+1}(\sinh w-\frac{1}{2}|z|^2
e^{-w})=0 \\
(n+1,j):\,k_{n+1j}+k_{n+1n+1}z_j e^{-w}=0\;,\,j=2,..,n\\
(n+1,n+1):\,-\sum_{j=2}^n k_{n+1j}z_j+k_{n+1n+1}(\cosh w-\frac{1}{2}|z|^2
e^{-w})=1 .$\\
It follows that $k_{n+1 n+1} e^{-w}=1$ and $z_j=-k_{n+1 j}$.
This determines $t$ and $x$. So for any $k\in K$ such that $k_{n+1
n+1}>0$ and any $a\in A,n\in N$ we can find $\tilde{k}\in\tilde{K},f\in F,
m\in N$ such that $k a n=\tilde{k} f m$.} $\Box$
\end{lem}

If $k_0$ is as above it is clear that $k_0\tilde{K} k_0=\tilde{K}$ so the set
$X:=\tilde{K}\cup k_0\tilde{K}$ is a (non connected) closed Lie subgroup of
$SO_0(1,n)$. Moreover from the lemma above any element of $SO_0(1,n)\ni g=k a
n$ such that $k_{n+1 n+1}\neq 0$ can be uniquely decomposed $g=x f m$ with
$x\in X,f\in F,m\in N$.\vspace{2ex}

\noindent{\bf Aknowledgments\hspace{1em}} I woud like to thank to Dr. S.
Zakrzewski for inspiration of this work and to Dr. A. Strasburger for
dicsussions. \\
The research was supported by Polish KBN grant No. 2 P301 020 07.

\ed